\documentclass[12pt,reqno]{amsart}
\usepackage{amsmath,amssymb,amsfonts,amscd,latexsym,amsthm,mathrsfs}
\usepackage[usenames]{color}
\usepackage[unicode]{hyperref}
\usepackage{graphicx}
\newcommand{\bbeta}{{\boldsymbol{\beta}}}
\newcommand{\bgamma}{{\boldsymbol{\gamma}}}

\begin{document}

\title{(Strange) gamma evaluations}

\author{Wadim Zudilin}
\address{Department of Mathematics, IMAPP, Radboud University, PO Box 9010, 6500~GL Nij\-me\-gen, The Netherlands}
\urladdr{https://www.math.ru.nl/~wzudilin/}

\date{}

\dedicatory{To Doron Zeilberger `The Creative' on his $3\times5^2$th birthday in year $(3^2\times5)^2$}


\begin{abstract}
We review `creative' strategies of closed-form evaluations of hypergeometric functions.
\end{abstract}

\maketitle

%==================================================

This exposition is a rather technical account of certain recent \emph{nice} developments in the field of `strange' or `gamma' evaluations of hypergeometric functions.
At the same time, we cover computational aspects neither deep nor systematic: the reader is invited to perform calculations alluded themself.

For details of the classical theory of hypergeometric functions and their $q$- (\emph{aka} `basic') analogues consult with \cite{Ba35,GR04,Sl66}. Examples of hypergeometric summations (and transformations)\,---\,those that happen to be gamma (but not strange) evaluations go back to at least the 18th century and are attached to the names as great as Euler and Gauss; the related binomial theorem and the Chu--Vandermonde summation are even earlier.
The adjective `strange' was coined after a list of Gosper's evaluations discussed by Gessel and Stanton at length in \cite{GS82} but, perhaps, first examples of such strange evaluations were systematically produced by Goursat in~\cite{Go81}. 

\section{Ebisu's methodology}
\label{sec1}

Some years ago, Ebisu \cite{Eb17} came up with a strategy to use contiguous relations of $_2F_1$ hypergeometric sums for proving numerous strange evaluations of the series.
Let us follow the explanation in \cite{BF22,Eb17}.

Consider a triple of hypergeometric parameters $a,b,c$ and abbreviate it by $\bbeta=(a,b,c)$;
denote
\[
F(\bbeta\mid z)={}_2F_1(a,b;c\mid z)
=\sum_{n=0}^\infty\frac{(a)_n(b)_n}{(c)_nn!}\,z^n
\]
the Euler--Gauss hypergeometric function \cite{Ba35,Sl66}, where $(a)_n=\Gamma(a+n)/\Gamma(a)$ is Pochhammer's symbol.
Let $k,l,m$ be a fixed triple of integers abbreviated $\bgamma=(k,l,m)$, the shift vector.
Using contiguity relations we can write down explicitly rational functions $R_\bgamma(\bbeta,z)$ and $Q_\bgamma(\bbeta,z)$ in $\mathbb{Q}(a,b,c,z)$ such that
\begin{equation}
\label{basic}
F(\bbeta+\bgamma\mid z)=R_\bgamma(\bbeta,z)F(\bbeta\mid z)+Q_\bgamma(\bbeta,z)F'(\bbeta\mid z),
\end{equation}
where $F'$ is the $z$-derivative of~$F$.
In other words, very shift of the parameters of a $_2F_1$ series by integers can be efficiently written as a linear combination of the original series and its $z$-derivative with rational-function coefficients.
The underlying algorithmic strategy makes use of composition of shifts of individual parameters by~1 listed in \cite[Sect.~1.4]{Sl66}.
A quadruple $(\bbeta,z_0)=(a,b,c,z_0)$ is called \emph{admissible} with respect to $\bgamma$
if $Q_\bgamma(\bbeta+t\bgamma,z_0)=0$ for all $t\in\mathbb C$. 
Choosing an admissible quadruple $(\bbeta,z_0)$ and substituting $\bbeta+t\bgamma$ for $\bbeta$ in~\eqref{basic} we obtain the functional equation
\begin{equation}\label{functional}
F(\bbeta+(t+1)\bgamma\mid z_0)=R_\bgamma(\bbeta+t\bgamma,z_0) F(\bbeta+t\bgamma\mid z_0)
\end{equation}
in parameter $t$. Viewing the coefficient $R_\bgamma(\bbeta+t\bgamma,z_0)$ as a rational function in~$t$, we factor it into linear terms:
\[
R_\bgamma(\bbeta+t\bgamma,z_0)=R_0\prod_{i=1}^r\frac{t+\alpha_i}{t+\delta_i},\quad R_0\in\mathbb C^*.
\]
Then
\[
R_0^t\prod_{i=1}^r\frac{\Gamma(t+\alpha_i)}{\Gamma(t+\delta_i)}
\]
satisfies the same functional equation as $F(\bbeta+t\bgamma\mid z_0)$.
It remains to identify these two functions of~$t$; this is a quite delicate part covered for most situations by \cite[Theorem~2.5]{BF22}. Note that the identification factor may also depend on~$t$.

An illustration given in \cite{BF22} deals with the case $(k,l,m)=(2,2,1)$ in which the vanishing of $Q_\bgamma(a+2t,b+2t,c+t,z)$ corresponds to that of
\begin{align*}
&
1 + a + b + a b - 2 c - a c - b c + c^2 + z + 2 a z + a^2 z + 2 b z +  a b z + b^2 z
\\ &\quad
- c z - a c z - b c z + (2 + a + b - 2 c + 7 z + 5 a z + 5 b z - 4 c z)t + (8z+1)t^2.
\end{align*}
The equations for the admissible quadruple are obtained by setting this polynomial in
$t$ identically zero; the corresponding system is solved by imposing
\[
z_0=-\tfrac18,\quad a=2t,\quad b=2t+\tfrac13,\quad c=t+\tfrac56
\]
or 
\[
z_0=-\tfrac18,\quad a=2t,\quad b=2t-\tfrac13,\quad c=t+\tfrac23.
\]
The first possibility leads to
\[
F\big(2(t+1),2(t+1)+\tfrac13,t+1+\tfrac56\mid -\tfrac18\big)
=\frac{16}{27}\times\frac{t+\tfrac56}{t+\tfrac23}\times F(2t,2t+\tfrac13,t+\tfrac56\mid -\tfrac18),
\]
while application of \cite[Corollary~2.8]{BF22} finally translates it into
\[
{}_2F_1\bigg(2t,2t+\frac13;t+\frac56 \Bigm| -\frac18\bigg)
=\bigg(\frac{16}{27}\bigg)^t\frac{\Gamma(t+\tfrac56)\Gamma(\tfrac23)}{\Gamma(t+\tfrac23)\Gamma(\tfrac56)}.
\]
A similar gamma evaluation is deduced for the second choice of admissible quadruple.

\section{Zeilberger's methodology}
\label{sec2}

There is another method to produce such closed-form evaluations.
It is similar in spirit and essentially equivalent to Ebisu's method discussed in Section~\ref{sec1}; it has been given by Zeilberger \cite{Ze05} a decade earlier.
A set of illustrative examples was provided by Shalosh B.~Ekhad in a separate piece \cite{Ek04}.
Zeilberger looks for a relation of the form
\begin{equation}\label{basic-Z}
F(\bbeta+(t+1)\bgamma\mid z)=\hat R_\bgamma(t,\bbeta,z)F(\bbeta+t\bgamma\mid z)+\hat Q_\bgamma(t,\bbeta,z)F(\bbeta+(t-1)\bgamma\mid z)
\end{equation}
instead of \eqref{functional}, then for conditions to eliminate the last term on the right-hand side and reduce the equality to
\begin{equation}\label{functional-Z}
F(\bbeta+(t+1)\bgamma\mid z_0)=\hat R_\bgamma(t,\bbeta,z_0) F(\bbeta+t\bgamma\mid z_0).
\end{equation}
Of course, at the end of the day the two functional relations \eqref{functional} and \eqref{functional-Z} tend out to be the same, so that solving for the vanishing of $\hat Q_\bgamma(t,\bbeta,z)$ is essentially the same task as doing this for $Q_\bgamma(\bbeta+t\bgamma,z)$.
There are some transparent differences though.

For Zeilberger's design \eqref{basic-Z} one does not need to know anything about contiguous relations: simply execute the Gosper--Zeilberger algorithm \cite{PWZ96} for the hypergeometric term
\begin{equation}\label{hyperterm}
A(t,n)=\frac{(a+kt)_n(b+lt)_n}{(c+mt)_n(1)_n}\,z^n
\end{equation}
to obtain the telescoping relation
\[
p_0(t)A(t+1,n)+p_1(t)A(t,n)+p_2(t)A(t-1,n)=B(t,n+1)-B(t,n)
\]
for some $B(t,n)$ which is a rational multiple of $A(t,n)$.
(The fact that the difference operator $p_0(t)T^2+p_1(t)T+p_2(t)$, with $T\colon t\mapsto t+1$ the shift operator, is of order at most~2 is guaranteed by the fact that $\sum_nA(t,n)$ is a $_2F_1$ hypergeometric series, in turn implied by the structure of contiguous relations discussed earlier.)
Summing over $n$ and checking the boundary behaviour of $B(t,n)$ at $n=0$ and $n=\infty$, leads to
\begin{equation}\label{basic-Z2}
p_0(t)F(\bbeta+(t+1)\bgamma\mid z)+p_1(t)F(\bbeta+t\bgamma\mid z)+p_2(t)F(\bbeta+(t-1)\bgamma\mid z)=0
\end{equation}
which is a form of \eqref{basic-Z}.
At the same time dealing with \eqref{basic-Z2} directly has its own advantages: for example, to reduce it to a shorter relation (that is, to a difference equation of order~1 as in \eqref{functional-Z}) we may look for either $p_0(t)=0$ or $p_2(t)=0$.
(These two strategies are essentially equivalent, as what comes out from executing one of them should be producible by the other.)
One can even look for $p_1(t)=0$; this leads to a first-order difference equation with respect to the shift $t\mapsto t+2$.
Furthermore, in most of the cases one can throw in an additional parameter $d$ into the hypergeometric term~\eqref{hyperterm}:
\begin{equation*}
%\label{hyperterm2}
A(t,n)=\frac{(a+kt)_n(b+lt)_n}{(c+mt)_n(d)_n}\,z^n
\end{equation*}
to still have the telescoping relation.
However, when $d\notin\mathbb Z_{>0}$, the difference equation in \eqref{basic-Z2} becomes in general inhomogeneous; solving it after reduction to order~1 requires extension of the results from~\cite{BF22}.

One more\,---\,purely technical\,---\,advantage of Zeilberger's scheme concerns finding appropriate $z_0$ for admissible quadruples. In both Ebisu's and Zeilberger's methods they come from the leading coefficients in the $t$ expansions of (the numerators of) the rational coefficients; in the form \eqref{basic-Z2} we just need the leading $t$-terms of $p_0(t)$ or $p_2(t)$ (in fact, of all three coefficients\,---\,this comes from a careful analysis of the hypergeometric structure).
The latter can be cheaply achieved for \emph{any} choice of $a,b,c$, so that we can first execute the algorithm of creative telescoping with the choice $a=b=c=0$ (well, we better take $c=1$ if $m$ coincides with either~$k$ or~$l$).
This leads to a shorter running time; after that we could run the algorithm again with $z_0$ specified but $a,b,c$ undetermined.

As a final remark, in spite of the obvious similarity of two methodologies, the webpage of \cite{Ek04} has the following addition on the 27th February 2014: ``Akihito Ebisu found many new amazing identities, using a novel new method'' and refers to what finally became \cite{Eb17}.

\section{WZ connections}
\label{sec3}

The fundamentals of the WZ machinery are well explained in the fundamental Wilf--Zeilberger paper \cite{WZ90}; the WZ acronym stands for the standard use of $w$ and $z$ as complex variables.
It is already made clear in \cite{Ze05} that an underlying mechanism behind gamma evaluations is via appropriate WZ pairs, that is, via hypergeometric terms $A(t,n)$ and $B(t,n)$ in two variables satisfying
\[
A(t+1,n)-A(t,n)=B(t,n+1)-B(t,n).
\]
In a general context of strange evaluations, Gessel \cite{Ge95} highlighted the importance of WZ-pairs already three decades ago.
This line has received a second birth in the recent work of Au \cite{Au24}, where he introduces a notion of WZ seed, with rules of generating infinitely many WZ pairs from a single one, and algorithmic strategies for realising a targetted identity in terms of a WZ seed.
The techniques from \cite{Au24} lead to proofs of numerous gamma evaluations which were open for years, in particular, of several Ramanujan-type formulas for (negative powers of) $\pi$;
these proofs in fact give multi-parameter extensions of such formulas for suitable choices of the differential parameter $z=z_0$, just like in Ebisu--Zeilberger evaluations but now for more general $_pF_{p-1}$ hypergeometric functions.
There is a limited supply for the choice of \emph{rational} $z_0$:%
\footnote{A mystery behind the choice of algebraic $z_0$ in the $_2F_1$ hypergeometric case is addressed in \cite[Remark 4.3]{BF22}:
such $z_0$ makes $x^l(1-x)^{m-l}(1-z_0x)^{-k}$ a \emph{Belyi map} in~$x$.
It remains open how this extends to the general hypergeometric situation.}
Gessel remarks in \cite[Sect.~9]{Ge95} about the 2-3-5 rule: one usually expects $z_0$ be admissible if both $z_0$ and $1-z_0$ involve the three smallest primes only in their factorisation (there are a few artificially designed exceptions though).
This is echoed in \cite{Au24} with a notice that even the prime $5$ does not originate from the WZ seeds coming from the classical hypergeometric summations.
Examples obtained in \cite{Ek04} or by the methods from \cite{Eb17,Ek04,Ze05} include gamma evaluations for
\[
{}_2F_1\big(-t,\tfrac12-t;\tfrac32+4t\mid\tfrac15\big), \quad
{}_2F_1\big(-t,\tfrac12-t;\tfrac52+4t\mid\tfrac15\big)
\]
and
\[
{}_2F_1\big(t,\tfrac12+t;6t\mid\tfrac45\big), \quad
{}_2F_1\big(t,\tfrac12+t;-1+6t\mid\tfrac45\big).
\]
Though they are weak as only depend on just one additional parameter~$t$, they can be extended to two-parameter (inhomogeneous!) gamma evaluations following a remark in Section~\ref{sec2}: each of the sums
\[
{}_3F_2\bigg(\begin{matrix} 1, \, -t, \, \tfrac12-t \\[2pt] d, \, 5d-\tfrac72+4t \end{matrix}\biggm|\frac15\bigg), \quad
{}_3F_2\bigg(\begin{matrix} 1, \, -t, \, \tfrac12-t \\[2pt] d, \, 5d-\tfrac52+4t \end{matrix}\biggm|\frac15\bigg)
\]
and
\[
{}_3F_2\bigg(\begin{matrix} 1, \, t, \, \tfrac12+t \\ d, \, 5-5d+6t \end{matrix}\biggm|\frac45\bigg), \quad
{}_3F_2\bigg(\begin{matrix} 1, \, t, \, \tfrac12+t \\ d, \, 4-5d+6t \end{matrix}\biggm|\frac45\bigg)
\]
solves a first order inhomogeneous difference equation (too lengthy though for displaying an example here).
Unlike the main results in \cite{BF22}, there is no general machinery to solve such inhomogeneous equations.
%So far no example with rational $z_0$ and corresponding to a gamma evaluation of $_2F_1$ or $_3F_2$ has been discovered.

\section{Higher rank reduction}
\label{sec4}

The two methods in \cite{Eb17} and \cite{Ze05} explore the same idea of forcedly reducing a rank~2 holonomic system to rank~1, using slightly different variants.
The work \cite{Eb17} furthermore outlines some possibilities for reducing higher-rank holonomic systems to rank~1, but not to rank 2 or higher (because of its design); this is discussed in \cite{Eb17} for $_3F_2$ (already known) evaluations and for Appell's $F_1$ double hypergeometric function.
Zeilberger's machinery can be extended to this setup along the lines.
In addition, one may source $_3F_2$ gamma evaluations from $_2F_1$ ones that match Clausen's formula
\[
{}_2F_1\big(a,b;a+b+\tfrac12\mid z\big)^2
={}_3F_2\bigg(\begin{matrix} 2a, \, 2b, \, a+b \\ 2a+2b, \, a+b+\tfrac12 \end{matrix}\biggm| z\bigg)
\]
or, more generally, Orr-type formulas \cite[Sect.~2.5]{Sl66}.
One example in \cite{BF22} corresponding to the shift $\bgamma=(-1,3,2)$ (see also \cite[p.~34, eq.~(ix)]{Eb17}) reads
\[
{}_2F_1\bigl(-t,3t+1,2t+\tfrac32\mid \tfrac14\bigr)
=\bigg(\frac{16}{27}\bigg)^t\frac{\Gamma(t+\frac54)\Gamma(t+\frac34)\Gamma(\frac76)\Gamma(\frac23)}{\Gamma(t+\frac76)\Gamma(t+\frac23)\Gamma(\frac54)\Gamma(\frac34)}
=\frac{2^{4t+2}\Gamma(2t+\frac32)\Gamma(\frac12)}{3^{3t+3/2}\Gamma(2t+\frac43)\Gamma(\frac23)};
\]
after Clausen-squaring it becomes
\begin{equation}
{}_3F_2\bigg(\begin{matrix} -2t, \, 6t+2, \, 2t+1 \\ 4t+2, \, 2t+\frac32 \end{matrix}\biggm| \frac14\bigg)
=\frac{2^{8t+4}\Gamma(2t+\frac32)^2\Gamma(\frac12)^2}{3^{6t+3}\Gamma(2t+\frac43)^2\Gamma(\frac23)^2}.
\label{Eb1}
\end{equation}
Similarly, the evaluation \cite[p.~36, eq.~(ix)]{Eb17}
\[
{}_2F_1\bigl(-t,3t+1;2t+\tfrac32\mid -\tfrac18\bigr)
=\frac{2^{3t+1}\Gamma(2t+\frac32)\Gamma(\frac12)}{3^{3t+1}\Gamma(t+\frac76)\Gamma(t+\frac56)}
\]
leads to
\begin{equation*}
{}_3F_2\bigg(\begin{matrix} -2t, \, 6t+2, \, 2t+1 \\ 4t+2, \, 2t+\frac32 \end{matrix}\biggm| -\frac18\bigg)
=\frac{2^{6t+2}\Gamma(2t+\frac32)^2\Gamma(\frac12)^2}{3^{6t+2}\Gamma(t+\frac76)^2\Gamma(t+\frac56)^2}.
%\label{Eb2}
\end{equation*}

Zeilberger's method \cite{Ze05} allows one to deal with a more general sum $a(t)=F(\bbeta+t\bgamma\mid z)$,
when $\bbeta$ and $\bgamma$ correspond to the $_pF_{p-1}$ hypergeometric function parameters.
The corresponding difference equation
\[
\sum_{i=0}^rp_i(t)a(t+i)=0
\]
of order $r\le p$ can be further reduced to a lower order when $p_r(t)=0$ (or $p_0(t)=0$) via an appropriate choice of the parameters. This does not (explicitly) rely on contiguous relations and give, for example, reductions of rank~3 holonomic systems to rank~2; the bottleneck is the running time.
Applying this method to some $_3F_2$-series (for which running time is finite) we found that
\[
{}_3F_2\bigg(\begin{matrix} \frac13, \, t, \, -\frac13+2t \\ \frac23, \, \frac12+t \end{matrix}\biggm| 2\bigg),
% (2^k Pochhammer[1/3, k] Pochhammer[t, k] Pochhammer[-(1/3) + 2 t, k])/(k! Pochhammer[2/3, k] Pochhammer[1/2 + t, k])
\quad
{}_3F_2\bigg(\begin{matrix} \frac23, \, t, \, -\frac23+2t \\ \frac43, \, \frac12+t \end{matrix}\biggm| 2\bigg),
% (2^k Pochhammer[2/3, k] Pochhammer[t, k] Pochhammer[-(2/3) + 2 t, k])/(k! Pochhammer[4/3, k] Pochhammer[1/2 + t, k])
\quad
{}_3F_2\bigg(\begin{matrix} \frac12, \, t, \, 1-t \\ 1, \, \frac12+t \end{matrix}\biggm| \frac14\bigg)
\]
satisfy second rather than third order difference equations in~$t$.
For example, if we denote the latter hypergeometric function by $G(t)$ then
\[
t(4t-1)(4t+1)G(t+1) - (2t+1)(10t^2-10t+3)G(t) + (t-1)(2t-1)(2t+1)G(t-1)=0.
\]
At $t=\frac12$, the series possesses the closed form
\begin{equation*}
{}_3F_2\bigg(\begin{matrix} \frac12, \, \frac12, \, \frac12 \\ 1, \, 1 \end{matrix}\biggm| \frac14\bigg)
=\frac{\sqrt3\,\Gamma(\frac13)^6}{2^{8/3}\pi^4}
%\label{Eb1s}
\end{equation*}
which is a special case, when $t=-\frac14$, of the evaluation~\eqref{Eb1}.

\section{$q$-Extensions}
\label{sec5}

There seems to be no straightforward way to extend the methodology from \cite{Ze05} and \cite{Eb17} to a $q$-setting;
as the examples below indicate the $q$-extensions rarely found for $_2F_1$ gamma evaluations come from a system of holonomic rank higher than~$2$, so their reduction to rank~$1$ is hardly explainable.
However such $q$-analogues of formulas from \cite{Eb17} do exist, and they are provable using classical basic hypergeometric techniques \cite{GR04} and creative telescoping \cite{PWZ96} or $q$-WZ pairs \cite{WZ90}.

Chern's identity \cite{Ch23}
\begin{equation}
\sum_{k=0}^\infty\frac{(q^2/a;q^2)_k(a;q)_k}{(q^3;q^3)_k(a^2q;q^2)_k}(-1)^kq^{\binom{k+1}2}a^k
=\frac{(aq;q^2)_\infty(a^3q^3;q^6)_\infty}{(a^2q;q^2)_\infty(q^3;q^6)_\infty}
\label{eq:id1}
\end{equation}
is one of such instances, which played a crucial role in the recent resolution of a combinatorial conjecture of Andrews and Uncu.
The identity can be proved by showing that the summand
\[
F_k(a,q)=\frac{(q^2/a;q^2)_k(a;q)_k}{(q^3;q^3)_k(a^2q;q^2)_k}(-1)^kq^{\binom{k+1}2}a^k
\]
satisfies the telescoping relation
\[
(1-a^2q)(1-a^6q^3)F_k(aq,q)-(1-a^4q)(1-a^4q^3)F_k(a,q)=G_{k+1}(a,q)-G_k(a,q)
\]
for a suitable choice of $G_k(a,q)$.

A similar strategy works for
\begin{equation}
\sum_{k=0}^\infty\frac{(q/a,a^2;q)_k}{(q^3,a^3q^2;q^3)_k}\,(1+aq^{2k})\,q^{k^2}a^k
=\frac{(-a;q)_\infty(a^2q;q^2)_\infty(q^2;q^3)_\infty}{(a^3q^2;q^3)_\infty}.
\label{eq:id2}
\end{equation}

Identity \eqref{eq:id1} is a $q$-analogue of Ebisu's formula \cite[p.~27, Sect.~4.3.3, eq.~(xiv)]{Eb17},
\[
\sum_{k=0}^\infty\frac{(1-a)_k(2a)_k}{k!\,(2a+\frac12)_k}\,\biggl(-\frac13\biggr)^k
=\frac{\Gamma(\frac12)\Gamma(2a+\frac12)}{3^a\Gamma(a+\frac12)^2},
\]
while identity \eqref{eq:id2} is a $q$-analogue of \cite[p.~27, eq.~(A.18)]{Eb17} (see also \cite[p.~696]{BF22}),
\[
\sum_{k=0}^\infty\frac{(1-a)_k(2a)_k}{k!\,(a+\frac23)_k}\,\biggl(\frac19\biggr)^k
=\frac{3^a\Gamma(\frac12)\Gamma(a+\frac23)}{4^a\Gamma(\frac23)\Gamma(a+\frac12)}.
\]

As L.~Wang points out, identity \eqref{eq:id1} is a specialisation of \cite[Exercise 3.29(i)]{GR04} (see also \cite[Corollary 6]{CW09}): write it as
\begin{align}
&
\sum_{k=0}^\infty\frac{(1-cq^{5k})(q^2/a,c;q^2)_k(a;q)_k(cq/a^2;q^3)_k(acq;q^6)_k}{(1-c)(q^3,acq;q^3)_k(a^2q;q^2)_k(cq^2/a,cq^4/a;q^4)_k}(-1)^kq^{\binom{k+1}2}a^k
\nonumber\\ &\quad
=\frac{(cq^2,aq;q^2)_\infty(cq^4/a^2,a^3q^3;q^6)_\infty}{(cq^2/a,a^2q;q^2)_\infty(q^3,acq^4;q^6)_\infty}
\label{eq:id3}
\end{align}
(notice the correction of typo in the product side) and choose $c=0$.
Now \cite[Exercises 3.29(ii), (iii)]{GR04} do not seem to possess suitable specialisations ($c=0$ is boring, and the series diverge as $a\to0$).
They suggest however the existence of a shifted version of \eqref{eq:id3} (which, unfortunately, does not follow from \cite[Proposition~5]{CW09}).

Specialisation $a=0$ of \cite[Corollary 10]{CW09} (or the shift $b\mapsto b/q^3$ in \cite[Corollary~9]{CW09}) reads
\begin{equation*}
\sum_{k=0}^\infty\frac{(q;q^2)_k(b;q^4)_k}{(q^4;q^4)_k(bq^3;q^6)_k}(-1)^kq^{k^2+2k}
=\frac{(b;q^4)_\infty(q^3;q^6)_\infty(q^{12};q^{12})_\infty}{(q^4;q^4)_\infty(bq^3;q^6)_\infty(b;q^{12})_\infty};
%\label{eq:id4}
\end{equation*}
this is a $q$-analogue of Ebisu's \cite[p.~26, eq.~(x)]{Eb17}.

Specialisation $a=0$ and $m\to\infty$ of \cite[Section~3.1, the first display]{CC21} leads to
\begin{equation*}
\sum_{k=0}^\infty\frac{(b,q/b;q)_k}{(q^2;q^2)_k(qd;q)_k}q^{\binom{k+1}2}d^k
=\frac{(qbd,q^2b/d;q^2)_\infty}{(qd;q)_\infty};
%\label{eq:id5}
\end{equation*}
this is a $q$-analogue of \cite[p.~22, eq.~(vii)]{Eb17}.

Specialisation $a=0$ of Rahman's \cite[eq.~(1.8)]{Ra93} results in
\begin{equation*}
\sum_{k=0}^\infty\frac{(b,d;q)_k}{(q;q)_k(qbd;q^2)_k}q^{\binom{k+1}2}
=\frac{(qb,qd;q^2)_\infty}{(q,qbd;q)_\infty}.
%\label{eq:id6}
\end{equation*}

Some other $q$-analogues of $_2F_1$ gamma evaluations can be found or extracted from the existing literature;
they however do not suggest a general method.
Perhaps, linking the corresponding $_2F_1$ hypergeometric sums with WZ pairs, connecting the latter with the classical WZ seeds from \cite{Au24} and performing the $q$-deformation of the latter as explained in Au's further work \cite{Au24b} will give one an algorithmic way of finding $q$-extensions.

\section{An attempt to conclude}
\label{sec6}

There are other techniques in use, somewhat algorithmic as well, to produce non-standard gamma evaluations of hypergeometric functions, for example, those discussed in \cite{CS13,GS82,Go81}.
Strange evaluations resulting from both Ebisu's and Zeilberger's methods automatically depend on an additional parameter~$t$. This limits their supply considerably. 
There is no evidence that the evaluations
\begin{align*}
{}_3F_2\biggl(\begin{matrix} \frac12, \, \frac12, \, \frac12 \\ 1, \, 1 \end{matrix} \biggm|\frac1{64}\biggr)
&=\frac{2}{7\pi}\prod_{j=1}^{7}\Gamma\biggl(\frac j7\biggr)^{\left(\frac j7\right)}
\\ \intertext{(this is the original example of CM-evaluation from \cite[\S\,8]{SC67}) and}
{}_3F_2\biggl(\begin{matrix} \frac12, \, \frac16, \, \frac56 \\ 1, \, 1 \end{matrix} \biggm|-\frac1{80^3}\biggr)
&=\frac{4\sqrt{15}}{43\pi}\prod_{j=1}^{42}\Gamma\biggl(\frac j{43}\biggr)^{\left(\frac j{43}\right)},
\end{align*}
where $\big(\frac{\cdot}{N}\big)$ denotes the Legendre symbol, possess $t$-parametric extensions.
Is there a method to prove such rigid identities?

\medskip
\noindent
\textbf{Acknowledgements.}
Discussions with many colleagues have crucially influenced my perception of gamma evaluations.
I am happy to thank Jean-Paul Allouche, Kam Cheong Au, Frits Beukers, Alin Bostan, Shane Chern, Vasily Golyshev, Jes\'us Guillera, Logan Kleinwaks, Christoph Koutschan, Christian Krattenthaler, Armin Straub, Ali Uncu and Doron Zeilberger for sharing interest with me in aspects of the topic.
I keep my memories about a mysterious 2019 talk on the method from \cite{Eb17} delivered by Akihito Ebisu at a conference in the T\^ohoku University, Sendai; it was not until an illuminating lecture by Frits Beukers on the results from his joint paper \cite{BF22} given online in 2021 when I could really absorb the ingredients and aesthetics of Ebisu's techniques.
Witnessing further developments of the WZ theory, also through $q$-glasses, convinced me to take another look into the subject.
I keep thinking on where all these things can go.

Finally, I am thankful to the Max-Planck Institute in Bonn (Germany) not only for financial support during my stay in April--June 2025 but also for a unique research atmosphere that helped me to crystallise my thoughts on the theme.

%==================================================


\begin{thebibliography}{99}

\bibitem{Au24}
\textsc{K.\,C. Au},
Wilf--Zeilberger seeds and non-trivial hypergeometric identities,
\emph{J. Symbolic Comput.} \textbf{130} (2025), Paper No.~102421, 34~pp.
%\emph{Preprint} (2023/24), 30~pp.; \href{https://arxiv.org/abs/2312.14051}{arXiv:\,2312.14051 [math.CO]}.

\bibitem{Au24b}
\textsc{K.\,C. Au},
Wilf--Zeilberger seeds: $q$-analogues,
\emph{Preprint} (2024), 28~pp.; \href{https://arxiv.org/abs/2403.04555}{arXiv:\,2403.04555 [math.CO]}.

\bibitem{Ba35}
\textsc{W.\,N.~Bailey},
\emph{Generalized hypergeometric series},
Cambridge Math. Tracts \textbf{32} (Cambridge Univ. Press, Cambridge, 1935).

\bibitem{BF22}
\textsc{F.~Beukers} and \textsc{J.~Forsg{\aa}rd},
$\Gamma$-evaluations of hypergeometric series,
\emph{Ramanujan J.} \textbf{58} (2022), 677--699.

\bibitem{CS13}
\textsc{M. Chamberland} and \textsc{A. Straub},
On gamma quotients and infinite products,
\emph{Adv. in Appl. Math.} \textbf{51} (2013), no. 5, 546--562.

\bibitem{CC21}
\textsc{X. Chen} and \textsc{W. Chu},
Hidden $q$-analogues of Ramanujan-like $\pi$-series,
\emph{Ramanujan J.} \textbf{54} (2021), no.~3, 625--648.

\bibitem{Ch23}
\textsc{S. Chern},
Asymmetric Rogers--Ramanujan type identities. I. The Andrews--Uncu conjecture,
\emph{Proc. Amer. Math. Soc.} \textbf{151} (2023), no.~8, 3269--3279.
%\href{https://doi.org/10.1090/proc/16332}{doi:10.1090/proc/16332}.

\bibitem{CW09}
\textsc{W. Chu} and \textsc{C.Y. Wang},
Partial sums of two quartic $q$-series,
\emph{SIGMA} \textbf{5} (2009), art.~050, 19~pp.

\bibitem{Eb17}
\textsc{A. Ebisu},
Special values of the hypergeometric series,
\emph{Mem. Amer. Math. Soc.} \textbf{248} (2017), no.~1177.

\bibitem{Ek04}
\textsc{S.\,B. Ekhad},
Forty ``strange'' computer-discovered [and computer-proved (of course!)] hypergeometric series evaluations,
\emph{Personal Journal of Shalosh B. Ekhad and Doron Zeilberger} (2004), 13~pp.;
\url{https://sites.math.rutgers.edu/~zeilberg/mamarim/mamarimhtml/strange.html}.

\bibitem{GR04}
\textsc{G.~Gasper} and \textsc{M.~Rahman},
\emph{Basic hypergeometric series},
2nd edition, Encyclopedia Math. Appl. \textbf{96}
(Cambridge Univ. Press, Cambridge, 2004).

\bibitem{Ge95}
\textsc{I.\,M. Gessel},
Finding identities with the WZ method,
%in: Symbolic computation in combinatorics $\Delta_1$ (Ithaca, NY, 1993)
\emph{J. Symbolic Comput.} \textbf{20} (1995), no.~5-6, 537--566.

\bibitem{GS82}
\textsc{I. Gessel} and \textsc{D. Stanton},
Strange evaluations of hypergeometric series,
\emph{SIAM J. Math. Anal.} \textbf{13} (1982), no.~2, 295--308.

\bibitem{Go81}
\textsc{É. Goursat},
Sur l'\'equation diff\'erentielle lin\'eaire qui admet pour int\'egrale la s\'erie hyperg\'eom\'etrique,
\emph{Ann. Sci. de l'\'ENS} (2) \textbf{10} (1881), 3--142.

\bibitem{PWZ96}
\textsc{M. Petkov\v sek}, \textsc{H.\,S. Wilf} and \textsc{D. Zeilberger},
$A=B$ (A K Peters, Ltd., Wellesley, MA, 1996).

\bibitem{Ra93}
\textsc{M.~Rahman},
Some quadratic and cubic summation formulas for basic hypergeometric series,
\emph{Canad. J. Math.} \textbf{45} (1993), no.~2, 394--411.

\bibitem{SC67}
\textsc{A. Selberg} and \textsc{S. Chowla},
On Epstein's zeta-function,
\emph{Crelle's J.} \textbf{227} (1967), 86--110.

\bibitem{Sl66}
\textsc{L.\,J. Slater},
\emph{Generalized hypergeometric functions}
(Cambridge University Press, Cambridge, 1966).

\bibitem{WZ90}
\textsc{D. Zeilberger} and \textsc{H.\,S. Wilf},
Rational functions certify combinatorial identities,
\emph{J. Amer. Math. Soc.} \textbf{3} (1990), no.~1, 147--158.

\bibitem{Ze05}
\textsc{D. Zeilberger},
Deconstructing the Zeilberger algorithm,
\emph{J. Difference Equations Appl.} \textbf{11} (2005), 851--856;
\url{https://sites.math.rutgers.edu/~zeilberg/mamarim/mamarimhtml/decon.html}.

\end{thebibliography}
\end{document}